\documentclass[a4paper,10pt]{article}

\usepackage[active]{srcltx}
\usepackage{bibgerm}
\usepackage[english]{babel}
\usepackage{tikz}
\usetikzlibrary{matrix}
\usetikzlibrary{arrows}
\usetikzlibrary{fit}
\usetikzlibrary{shapes}

\usepackage[T1]{fontenc}
\usepackage{amssymb}
\usepackage[latin1]{inputenc}
\usepackage{amsmath}
\usepackage{amsthm}
\usepackage{amsfonts}
\usepackage{a4wide}

\newtheorem{defi}{Definition}[section]
\newtheorem{remark}[defi]{Remark}
\newtheorem{example}[defi]{Example}
\newtheorem{theorem}[defi]{Theorem}
\newtheorem{lemma}[defi]{Lemma}
\newtheorem{corollary}[defi]{Corollary}

\renewcommand{\leq}{\leqslant}
\renewcommand{\geq}{\geqslant}
\newcommand{\modC}{\operatorname{\mathbf{mod}}}

\newcommand{\length}{\operatorname{length}}

\newcommand{\Hom}{\operatorname{Hom}}
\newcommand{\Ext}{\operatorname{Ext}}
\newcommand{\Jac}{\operatorname{Jac}}

\newcommand{\sign}{\operatorname{sgn}}
\newcommand{\Tr}{\operatorname{Tr}}

\newcommand{\Ker}{\operatorname{Ker}}

\newcommand{\Irr}{\operatorname{Irr}}
\newcommand{\Rad}{\operatorname{Rad}}
\newcommand{\Soc}{\operatorname{Soc}}

\newcommand{\id}{\operatorname{id}}

\newcommand{\eps}{\varepsilon}

\newcommand{\Sym}{\Sigma} 

\newcommand{\iso}{\cong}

\newcommand{\OO}{\mathcal{O}}

\newcommand{\Q}{\mathbb{Q}}
\newcommand{\Z}{\mathbb{Z}}
\newcommand{\F}{\mathbb{F}}
\newcommand{\N}{\mathbb{N}}

\usepackage{ctable}

\begin{document}

\title{Basic Orders for Defect Two Blocks of $\Z_p\Sym_n$}
\author{Florian Eisele\\ \\ \small\textit{Lehrstuhl D f\"ur Mathematik, RWTH Aachen, Templergraben 64, 52062 Aachen, Germany}\\\small\texttt{florian.eisele@rwth-aachen.de}}
\date{}

\maketitle

\begin{abstract}
	We show how basic orders for defect two blocks of symmetric groups over the ring of $p$-adic integers can
	be constructed by purely combinatorial means.
\end{abstract}

\section{Introduction and Outline}

Blocks of group rings of symmetric groups of small defect have been
subject to extensive study in the past. In this paper we are concerned with blocks of defect two.
For those, it has been shown in \cite{ScopesDef2}  
that (among other things) the decomposition numbers are all $\leq 1$. This bound for the decomposition numbers
is of particular interest in the context of the theory developed in \cite{Plesken}. 
Namely,
if $B$ is a defect two block of $\Z_p\Sym_n$, then
the image of $B$ under an irreducible representation of $\Q_p \otimes_{\Z_p}B$ is a so-called
graduated order (also known as a tiled order or a split order). These can be described easily in terms of only a few numerical invariants. 
Therefore describing these images is what we do first. This yields an overorder $\Gamma$ of a basic order $B_0$ of
$B$. We then go on to determine how $B_0$ is embedded in $\Gamma$. Here the shape of the $\Ext$-quiver and the decomposition matrix play an important role,
as these will turn out to control to what extent we can modify generators of $B_0$ by conjugation.
We find that those considerations determine $B_0$ uniquely up to conjugacy. $B_0$ is then given by generators 
in a direct sum of matrix rings over $\Q_p$.
It should not be too hard, in any particular case, to derive from this a presentation as a quiver algebra of $\F_p \otimes_{\Z_p}B_0$.

For the principal block of $\Z_p\Sym_{2p}$, the basic orders have been determined  in
\cite{NebeZpS2p}, and our approach is a generalization of that. The aforementioned paper relies, however,
heavily on the explicit knowledge of (among other things) the decomposition matrix of the blocks treated in it.
We, on the other hand, get by without such explicit information, which allows us to treat \emph{all} defect two blocks 
of symmetric groups (and makes it possible for the reader to verify the proofs without inspecting any large tables). In fact,
the following is a list of the known properties of defect two blocks of symmetric groups that we are going to use:
\begin{remark}[Known facts]\label{knownfacts}
	\begin{enumerate}
	 \item[(i)] The decomposition numbers of a defect two block of a symmetric group are all $\leq 1$
		and the off-diagonal Cartan numbers are all $\leq 2$.
		The diagonal Cartan numbers are all $\geq 3$. (see \cite{ScopesDef2}).
	\item[(ii)] The decomposition matrix of a defect two block of a symmetric group can be computed by combinatorial means, 
		for instance using the Jantzen-Schaper formula (see \cite{BensonDecomp}).
	\item[(iii)] 
	For any two simple modules $S$ and $T$ in a defect two block of $k\Sym_n$ we have
	$\dim_k \Ext^1_{k\Sym_n}(S,T)=\dim_k \Ext^1_{k\Sym_n}(T,S) \leq 1$ (see \cite{ScopesDef2}).
	We will in particular consider $\Ext$-quivers of defect two blocks of $k\Sym_n$ as \emph{undirected} graphs, as for each edge, there is an edge going back. We will
	represent those two directed edges  by a single \emph{undirected} edge.
	\item[(iv)] The $\Ext$-quiver of a defect two block of a symmetric group is a bipartite graph according to \cite[Theorem 3.2]{ChuangTanYoung}. Bipartite means,
		in this context, that the set of vertices of the quiver can be partitioned in two parts such that any edge connects vertices coming from
		different parts.
	\end{enumerate}
\end{remark}

\section{Notation}

Let $p$ be a prime and $(k,\OO,K)$ be a $p$-modular system such that $K$ is unramified over $\Q_p$.
For any $n \in \N$ we denote by $\Sym_n$ the symmetric group on $n$ points. 

In this paper, all modules are, unless stated otherwise, right modules. Whenever $\Lambda$ is an $\OO$-order, 
$V$ is a $K\otimes_\OO \Lambda$-module and $S$ is a simple $k\otimes_\OO \Lambda$-module, we denote by 
$[V:S]$ the multiplicity of $S$ in $k\otimes_\OO L$, where $L$ is any full $\Lambda$-lattice in $V$. If $V$ is simple as well, then $[V:S]$ is just the decomposition number associated to $V$ and $S$.

Whenever we have an $\OO$-lattice $L$ in a $K$-vector space $V$ that carries a symmetric bilinear form
$T: V\times V \rightarrow K$, we define its \emph{dual} $L^\sharp$ to be the $\OO$-lattice $\{l\in V\ | \ T(l,L) \subseteq \OO \}$.
We say $L$ is self-dual if $L=L^\sharp$.

Our notation concerning the representation theory of symmetric groups is essentially as in \cite{James}.
In particular, for a partition $\lambda$ 
of $n$ and a commutative ring $R$ we denote by $S^\lambda_R = S^\lambda_\Z \otimes_\Z R$ the corresponding Specht-module. 
If $\mu$ is a $p$-regular partition of $n$, we denote by $D^\mu$ the corresponding simple module defined over $k$,
that is, $D^\mu = S^\mu_k / \Rad S^\mu_k$. 

For a partition $\lambda$ of $n$, we denote by $\lambda^\top$ its transposed. For a $p$-regular partition $\mu$ of $n$
we denote by $\mu^M$ its image under the Mullineux-map, that is, define $\mu^M$ so that $D^{\mu^M} \iso D^\mu \otimes_k \sign$ holds.

In addition to that, we will use the following (non-standard) notation:

\begin{defi}
 	Define for a $p$-regular
	partition $\mu$ of  $n$ the set
	\begin{equation}
 		c_\mu := \left\{ \lambda \textrm{ a partition of $n$}\ \mid\ [S^\lambda_K : D^\mu] \neq 0  \right\}
	\end{equation}
	Define for any partition $\lambda$ of $n$ the set
	\begin{equation}
		r_\lambda := \left\{ \mu \textrm{ a $p$-regular partition of $n$}\ \mid\ [S^\lambda_K : D^\mu] \neq 0 \right\}
	\end{equation}
\end{defi}

\section{Some Facts about Orders}

In this section we recollect some facts about $\OO$-orders, mostly about graduated orders as
defined in \cite{Plesken}. All statements made in this section can be found either in \cite{Plesken} or in \cite{NebeH}.
We specialize everything to the splitting field case, as that case is what we need 
for the symmetric group blocks later on.

\begin{defi}
	Let $A$ be a semisimple $K$-split $K$-algebra and $\Lambda \subset A$ be a full $\OO$-order in $A$. Then $\Lambda$ is called 
	\emph{graduated} if $\Lambda$ contains a full set $e_1,\ldots,e_n$ of orthogonal idempotents which 
	are primitive in $A$.
\end{defi}

\begin{remark}
	Clearly a graduated order $\Lambda$ also contains the central primitive idempotents $\eps_1,\ldots,\eps_h$ of $A$.
	Thus we have \begin{equation}
             \Lambda = \bigoplus_i \eps_i\Lambda
            \end{equation}
	Each $\eps_i \Lambda$ is by itself a graduated order in $\eps_iA$, which means
	that it suffices to describe graduated orders in simple $K$-algebras.
\end{remark}

\begin{theorem}
	Any graduated order $\Lambda \subset K^{n\times n}$ (for some $n\in \N$)
	is conjugate to an order of the form
	\begin{equation}
		\Lambda(\OO, \hat{m}) = \bigoplus_{i,j=1}^n \left\langle p^{\hat{m}_{i,j}} \cdot e_{i,j}\right\rangle_\OO
	\end{equation}
	where $\hat{m} \in \Z_{\geq 0}^{n\times n}$ is some matrix, and $e_{i,j}$ denotes the element in $K^{n\times n}$ that 
	has a ``$1$'' at the $(i,j)$-entry and zeros elsewhere.
	
	By a further conjugation with a permutation matrix we may moreover achieve that all of the following holds:
	There is an integer $v$, a vector $d\in \Z_{>0}^v$ and a matrix $m\in \Z_{\geq 0}^{v\times v}$ such that
	$\sum_{i=1}^v d_i = n$ and 
	\begin{equation}
		\hat{m}_{i,j} = m_{s,t} \quad \forall i,j \textrm{ with } \sum_{l=1}^{s-1} d_l < i \leq \sum_{l=1}^s d_l
		\textrm{ and } \sum_{l=1}^{t-1} d_l < j \leq \sum_{l=1}^t d_l
	\end{equation}
	where the matrix $m$ satisfies
	\begin{enumerate}
	 \item[(i)] $m_{i,j} + m_{j,k} \geq m_{i,k} \quad \forall i,j,k \in \{1,\ldots,v\}$
	 \item[(ii)] $m_{i,j}+m_{j,i} > 0 \quad \forall i,j \in \{1,\ldots,v\}$
	\end{enumerate}
	We call the matrix $m$ an \emph{exponent matrix} for $\Lambda$. We call the vector $d$ a \emph{dimension
	vector}.
\end{theorem}

\begin{defi}
	Any matrix $m\in\Z_{\geq 0}^{v\times v}$ subject to the two conditions above, together with a vector
	$d\in\Z_{>0}^v$, defines a graduated order, which we denote by $\Lambda(\OO,m,d)$.
\end{defi}

\begin{remark}
	The integer $v$ from above is precisely the number of isomorphism classes of simple $\Lambda$-modules.
	The entries of the dimension vector $d$ are equal to the $k$-dimensions of those simple
	modules. In particular it is easily seen that $\Lambda(\OO,m,(1,\ldots,1))$ is a basic order of 
	$\Lambda(\OO,m,d)$.
\end{remark}

\begin{remark}\label{remark_simple_graduated}
	By fixing an exponent matrix we in particular fix a bijection
	\begin{equation}
		\{ 1,\ldots, v\} \leftrightarrow \{\textrm{ Isomorphism classes of simple $\Lambda$-modules }\}
	\end{equation}
	This bijection is afforded by the following map:
	\begin{equation}
		i \mapsto e_{s,s} \Lambda / \Rad e_{s,s}\Lambda \quad \textrm{with $s$ subject to } \sum_{l=1}^{s-1} d_l < s \leq \sum_{l=1}^{s}d_l
	\end{equation}
	Note furthermore that the projective cover of a simple module is an irreducible lattice. Therefore
	said projective cover is the \emph{unique} lattice with top isomorphic to that particular simple module. 
\end{remark}

The next theorem is the main reason why we are interested in graduated orders.
\begin{theorem}\label{thm_decomp_01}
	Let $A$ be finite-dimensional, semisimple $K$-split $K$-algebra. Let $\Lambda\subset A$ be a full $\OO$-order
	such that $k$ splits $k\otimes_\OO\Lambda$.
	If $V$ is a simple $A$-module and $\eps\in Z(A)$ is the corresponding central primitive idempotent, then
	$\eps\Lambda$ is a graduated order if and only if the decomposition numbers $[V : S]$ are 
	$\leq 1$ for all simple $\Lambda$-modules $S$ (i. e., $p$-reductions of irreducible lattices are multiplicity-free). 
\end{theorem}

\begin{remark}\label{remark_exp_mat_selfdual}
	Let $G$ be a finite group such that $k$ and $K$ are both splitting fields for $G$, let $\chi \in \Irr_K(G)$ be an irreducible character with decomposition numbers 
	$\leq 1$ and let $\eps_\chi \in Z(KG)$ be the
	corresponding central primitive idempotent. As seen above we have $\eps^\lambda \OO G \iso \Lambda(\OO,m,d)$ for some exponent matrix $m$ and dimension vector $d$. 
	It can be shown that $\OO G$ being a self-dual order with respect to the regular trace implies the following inequality:
	\begin{equation}
		m_{i,j}+m_{j,i} \leq \nu_p\left( \frac{\chi(1)}{|G|} \right) \quad \forall i,j
	\end{equation}
\end{remark}

\begin{theorem}[{\cite[Corollary 24]{NebeAbramenko}}]\label{thm_invol_idempot}
	Let $\Lambda$ be an $\OO$-order that carries an involution (i. e., an anti-automorphism of order two) $^\circ: \Lambda \rightarrow \Lambda$. 
	Assume moreover that $2\in\OO^\times$. Then 
	there is a full set of primitive pairwise orthogonal idempotents $e_1,\ldots,e_n \in \Lambda$
	and an involution $\sigma\in \Sym_n$ such that $e_i^\circ = e_{\sigma(i)}$.
\end{theorem}

\begin{corollary}\label{corollary_involution_basic_algebra}
	If $\Lambda$ is as in Theorem \ref{thm_invol_idempot}, then one may choose an idempotent $e\in \Lambda$ with
	$e^\circ = e$ such that $e\Lambda e$ is a basic algebra for $\Lambda$. In particular, a basic algebra for
	$\Lambda$ may be assumed to carry an involution as well. It should be noted that the assumption ``$2\in \OO^\times$'' of the preceding theorem is not
	necessary for this corollary to hold. 
\end{corollary}

\begin{theorem}\label{thm_sconst_invol}
	If $\Lambda = \Lambda(\OO,m,d)$ is a graduated order that carries an involution $^\circ: \Lambda \rightarrow \Lambda$ then 
	there is an involution $\sigma \in \Sym_v$ such that
	\begin{equation}
		m_{i,j} + m_{j,k} - m_{i,k} = m_{\sigma(k),\sigma(j)} + m_{\sigma(j),\sigma(i)} - m_{\sigma(k),\sigma(i)} \quad\forall i,j,k\in\{1,\ldots,v\}
	\end{equation}
	and $d_i=d_{\sigma(i)}$ for all $i\in\{1,\ldots,v\}$.
\end{theorem}

\begin{remark}\label{remark_involution}
	An involution $-^\circ$ on an order $\Lambda$ clearly induces an equivalence between the categories 
	$\modC_\Lambda$ and ${_\Lambda}\modC$. By slight abuse of notation, we denote this equivalence 
	(and its inverse) by $-^\circ$ as well. Then the functor $M \mapsto \Hom_\OO(M,\OO)^\circ$ is 
	an auto-equivalence of $\modC_\Lambda$. It hence permutes the simple $\Lambda$-modules.
	In the situation of the last theorem, in combination with Remark \ref{remark_simple_graduated}, 
	 $M \mapsto \Hom_\OO(M,\OO)^\circ$ does in particular induce a permutation on the set
	$\{1,\ldots,v\}$. The involution $\sigma$ in the Theorems \ref{thm_invol_idempot} and \ref{thm_sconst_invol} may 
	be chosen to equal this permutation.
	In particular, if $\Lambda = \eps^\lambda \OO\Sym_n$ (for some partition $\lambda$ of $n$) 
	is a graduated order, then we may choose $\sigma = \id$.
\end{remark}

\section{A Theorem on Amalgamation Depths}

The following theorem will be used only in a very special case. 
It is useful in many situations though,
so we give the general version.
In particular, we drop the assumption that $K$ be unramified over $\Q_p$, and denote by $\pi$
a generator for the maximal ideal of $\OO$.
 It is a fact about symmetric orders in $K^n$, taken as an algebra with component-wise multiplication (i. e. a semisimple algebra with $n$ non-isomorphic simple modules of dimension one). This setup is in a way orthogonal to what we looked at above, where we concentrated on orders
in $K^{n\times n}$ (i. e. a semisimple algebra with a single simple module of dimension $n$).
\begin{theorem}\label{thm_amal_comm}
	Let $\Lambda$ be a local symmetric suborder of the commutative $\OO$-order $\OO^n$. Fix a  $K^n$-equivariant symmetric bilinear form $K^n \times K^n \rightarrow K$  such that $\Lambda = \Lambda^\sharp$ with respect to that form. By $\eps_i$ we denote the $i$-th standard basis vector in $K^n$. Then we claim:
	If $L \leq K^n$ is a full $\Lambda$-lattice with 
	\begin{equation}\label{eqn_amal_comp_free}
\frac{L\cdot\eps_i}{L\cdot\eps_i\cap L} \iso_\OO 
	\frac{\Lambda\cdot\eps_i}{\Lambda\cdot\eps_i\cap \Lambda} \quad\textrm{for \emph{some} } i \in \{1,\ldots, n\}	 
	\end{equation}
	then $L \iso_\Lambda \Lambda$.
	\begin{proof}
		Let $L$ be a full $\Lambda$-lattice in $K^n$ not isomorphic to $\Lambda$. 
		Without loss we may assume that $L \subseteq \OO^n$ and $L\cdot \eps_i = \OO \cdot \eps_i$ for all $i\in\{1,\ldots,n\}$. We are going to  show that there can be no counter-example to the statement of the theorem, that is, that for each possible $L\ncong \Lambda_\Lambda$ we have that (\ref{eqn_amal_comp_free}) does not 
		hold for any $i\in\{1,\ldots,n\}$.

		First assume $\Lambda \subsetneqq L \subseteq \OO^n$. 
	 	We then have $L^\sharp \subseteq \Jac(\Lambda)$, and hence $L \supseteq \Jac(\Lambda)^\sharp$. Now, since $\Lambda$ is local and symmetric, the following holds:
		\begin{equation}
			\frac{\Jac(\Lambda)^\sharp\cdot \pi}{\Lambda\cdot \pi} = \Soc\left( \frac{\Lambda}{\Lambda\cdot\pi} \right)
		\end{equation}
		Therefore: If $l \in \Lambda$ such that $l + \Lambda\cdot \pi \in \Soc(\Lambda/\Lambda\cdot\pi)$
		then $l\cdot \pi^{-1} \in L$.

		Now let $l\in \Lambda\cdot\eps_i \cap\Lambda$ (where $i$ is arbitrary) such that $l\notin \Lambda\cdot\pi$. Then
		$\left(l\cdot \Lambda + \Lambda\cdot\pi \right) / \Lambda\cdot\pi \iso_\Lambda k$ (where $k$ is viewed as the simple $\Lambda$-module). This implies that $l + \Lambda\cdot \pi \in \Soc(\Lambda/\Lambda\cdot\pi)$, and 
		thus according to the above $l\cdot \pi^{-1} \in L$. Since $L\cdot\eps_i = \Lambda \cdot\eps_i = \OO \cdot \eps_i$, we conclude
		\begin{equation}
			\length_\OO \frac{L\cdot\eps_i}{L\cdot\eps_i\cap L} \leq \length_\OO \frac{\Lambda\cdot\eps_i}{\Lambda\cdot\eps_i\cap \Lambda} - 1 
		\end{equation}
		and this holds for each $i$, as $i$ was chosen arbitrary.
		One implication of the above is that for any idempotent $1 \neq \eps \in K^n$
		and any $i$ with $\eps \cdot \eps_i \neq 0$ the epimorphism 
		\begin{equation}\label{proper_epi}
			\Lambda\cdot\eps_i/\Lambda\cdot\eps_i\cap\Lambda \twoheadrightarrow \Lambda\cdot\eps_i / (\Lambda\cdot \eps_i) \cap (\Lambda\cdot\eps)
		\end{equation} 
		is proper. 
		Also we have shown at this point that a $\Lambda$-lattice $L$ with $\Lambda \subseteq L \subseteq \OO^n$ cannot possibly be a counterexample to the statement of the theorem.

		Now we consider an arbitrary $\Lambda$-lattice $L\subseteq \OO^n$ with $L\cdot\eps_i = \OO\cdot\eps_i$ for all $i$. We pick an element $v\in L$ with $v\cdot \eps_1 = \eps_1$. If $v\cdot \eps_i \in \OO^\times\cdot \eps_i$ for all $i$ then clearly $\Lambda \subseteq v^{-1}\cdot L \subseteq \OO^n$, and what we have shown 
		above implies that  $L$ is not a counterexample to the theorem.
		So assume that there is a $j\in\{2,\ldots,n\}$ such that $v\cdot\eps_j = r\cdot \eps_j$ with $r\in (\pi)_\OO$. Then we pick a $w\in L$ with $w\cdot \eps_j = \eps_j$ and look at $v' = v-r\cdot w$.
		By construction $v'\cdot\eps_1 \in \OO^\times \cdot \eps_1$ and $v'\cdot \eps_j = 0$.
		Now we have a series of epimorphisms
		\begin{equation}
			\frac{\Lambda\cdot\eps_1}{\Lambda\cdot\eps_1 \cap \Lambda} \twoheadrightarrow 
			\frac{\Lambda\cdot\eps_1}{\Lambda\cdot\eps_1 \cap \Lambda\cdot(1-\eps_j)} \twoheadrightarrow 
			\frac{L\cdot\eps_1}{L\cdot\eps_1\cap\Lambda\cdot v'} \twoheadrightarrow \frac{L\cdot \eps_1}{L\cdot\eps_1\cap L}
		\end{equation}
		of which at least the first one is proper, since it is a special case of the epimorphism in (\ref{proper_epi}). Hence the leftmost and the rightmost term cannot possibly be
		isomorphic. Repetition of this argument with $\eps_1$ replaced by $\eps_2,\eps_3,\ldots,\eps_n$
		yields that the theorem holds for $L$.
	\end{proof}
\end{theorem}

\begin{remark}[Applications to our situation]\label{remark_amal_comm}
	Let $\Lambda$ be an $\OO$-order in a semisimple $K$-split $K$-algebra $A$. Let $e$ be a primitive idempotent in $\Lambda$. Assume moreover that
	the decomposition numbers of $\Lambda$ are all $\leq 1$. Then $e\Lambda e$ is isomorphic to a local $\OO$-order in some $K^n$.
	Let $\eps_1,\ldots,\eps_h$ be the central primitive idempotents in $A$, and assume that $\Lambda=\Lambda^\sharp$ with respect
	to the trace bilinear form
	\begin{equation} 
		T_u: A\times A \rightarrow K: \ (a,b) \mapsto \sum_{i=1}^h\Tr(\eps_i\cdot u_i\cdot a\cdot b) \quad \textrm{ for some elements } u_i \in K \setminus (\pi)_\OO
	\end{equation}
	Then, by elementary linear algebra, we have
	\begin{equation}
		\frac{\eps_i \cdot e\Lambda e}{\eps_i \cdot e\Lambda e \cap e\Lambda e} \iso_\OO \OO/u_i^{-1}\OO 
	\end{equation}
	Let $I \subseteq \{1,\ldots,h\}$ be some set of indices such that $e\cdot\eps_i\neq 0 \ \forall i\in I$ and
	$e\cdot \eps \neq e$, where we put $ \eps := \sum_{i\in I} \eps_i$. Then $\eps\cdot e \Lambda e \oplus (1-\eps)\cdot e\Lambda e \ncong_{e\Lambda e} e \Lambda e$.
	Hence 
	\begin{equation}
		\length_\OO \frac{\eps_i\cdot e\Lambda e}{(\eps_i\cdot e\Lambda e)\cap (\eps \cdot e\Lambda e)} \leq \nu_\pi(u_i^{-1})-1\quad\forall i\in I	 
	\end{equation}
	Now we specialize to the (unramified) defect two case. This case corresponds to $\nu_\pi(u_i)=-2$ for all $i\in I$. We hence have 
	\begin{equation}\label{am1algebra}
		\eps \cdot e\Lambda e \iso \Gamma_{|I|} := \langle (1,1,1,\ldots,1), (0,\pi,0,\ldots,0),(0,0,\pi,\ldots,0),(0,0,0,\ldots,\pi) \rangle_\OO \subset K^{|I|}
	\end{equation}
	Formula (\ref{am1algebra}) is the consequence of the last theorem that we are actually going to use later.
	In the same vein is the following: One easily concludes from the last theorem that the radical idealizer of $e\Lambda e$ (that is, the largest subset $\Gamma \subset K^n$ such that $\Gamma \cdot \Jac(e\Lambda e) \subseteq  \Jac(e\Lambda e)$), which 
	for self-dual orders is known to be equal to $\Jac(e\Lambda e)^\sharp$ (see \cite{NebeRadicalId}), is isomorphic to the algebra $\Gamma_{|J|}$ as defined in (\ref{am1algebra}) with 
	$J=\{i\ |\ \eps_i\cdot e \neq 0\}$. Thus 
	\begin{equation}
		e\Lambda e \cong \langle (1,\ldots,1)\rangle_\OO + \Gamma_{|J|}^\sharp
	\end{equation}
	which basically says that $e\Lambda e$ is already determined by the $u_i$ (which in case $\Lambda$ is a group ring are just 
	the character degrees divided by the group order). Of course these defect two results have more elementary proofs.
\end{remark}

\section{Defect Two Blocks of Symmetric Groups}

\begin{remark}
 	In what follows, $p$ will always be an odd prime. When we say that a partition is in a defect two block, that simply means that
	it is of $p$-weight two.
\end{remark}

\begin{defi}[Jantzen-Schaper-Filtration]
	Let $\lambda$ be a partition of some $n\in\N$, and let
	\begin{equation}
		(-,=):\ S^\lambda_{\OO} \times S^\lambda_\OO \rightarrow \OO
	\end{equation}
	be the natural bilinear form on $S^\lambda_\OO$ inherited from the permutation module $M^\lambda_\OO$
	(see \cite{James} for details). Then we define for $i\in \Z_{\geq 0}$
	\begin{equation}
		S^\lambda_\OO(i) := \left\{ m \in S^\lambda_\OO \ | \ (m, S^\lambda_\OO) \subseteq p^i\cdot\OO \right\}
	\end{equation}
	and 
	\begin{equation}
		S^\lambda_k(i) := \frac{S^\lambda_\OO(i) +  p \cdot S^\lambda_\OO}{ p \cdot S^\lambda_\OO} \leq S^\lambda_k
	\end{equation}
	The filtration $S^\lambda_k = S^\lambda_k(0) \geq S^\lambda_k(1) \geq S^\lambda_k(2) \geq \ldots$
	is called the \emph{Jantzen-Schaper filtration} of $S^\lambda_k$.
\end{defi}

\begin{remark}\label{remark_js_semisimple}
	If $S^\lambda_k$ is multiplicity-free, then all layers $S^\lambda_k(i)/S^\lambda_k(i+1)$ of the Jantzen-Schaper
	filtration are semisimple. So, in particular, this holds for an $S^\lambda_k$ in a defect two block.
	\begin{proof}
		Consider the restriction of the standard bilinear form $(-,=)$ on $S^\lambda_\OO$ to $S^\lambda_\OO(i)$ for 
		some $i$. By definition of $S^\lambda_\OO(i)$, this takes values in $(p^i)_\OO$. Thus we may look at
		$p^{-i}\cdot (-,=)$, which defines a bilinear form on $S^\lambda_\OO(i)$ with values in $\OO$.
		We reduce this modulo $p$ to get a bilinear form on $k\otimes_\OO S^\lambda_\OO(i)$. 
		Clearly 
		\begin{equation}X := \frac{k\otimes_\OO S^\lambda_\OO(i)}{k\otimes_\OO S^\lambda_\OO(i)\cap k\otimes_\OO S^\lambda_\OO(i)^\perp}
		\end{equation}
 is a self-dual $k\Sym_n$-module. Since $S^\lambda_k$ (and therefore
		also $k\otimes_\OO S^\lambda_\OO(i)$) is multiplicity-free, and all simple $k\Sym_n$-modules
		are self-dual, we must hence have that $X$ is semisimple (as any simple module occurring in the radical
		would otherwise turn up again in the socle, giving it a multiplicity of at least two).
		Now we have the natural epimorphism $S^\lambda_\OO(i) \twoheadrightarrow S^\lambda_k(i)$, 
		giving rise to an epimorphism $k\otimes_\OO S^\lambda_\OO(i) \twoheadrightarrow S^\lambda_k(i)$.
		This epimorphism maps $k\otimes_\OO S^\lambda_\OO(i)\cap k\otimes_\OO S^\lambda_\OO(i)^\perp$ into $S^\lambda_k(i+1)$, which is best seen by diagonalizing the bilinear form on $S^\lambda_\OO$.
		Thus we get an epimorphism $X \twoheadrightarrow S^\lambda_k(i)/S^\lambda_k(i+1)$, implying that 
		the latter is also semisimple.
	\end{proof}
\end{remark}

The following theorem essentially summarizes what can be said about the structure of Specht modules in defect two blocks.
\begin{theorem}\label{thm_structure_schaper}
	Let $\lambda$ be a partition in a defect two block. Then the 
	Jantzen-Schaper-quotients of $S^\lambda_k$ and $S^{\lambda^\top}_k$ may be described as follows:
	\begin{enumerate}
	 \item[(i)] If $\lambda$ and $\lambda^\top$ are both $p$-regular, then
		$$
			\begin{array}{rclrcl}
			S^\lambda_k(0)/S^\lambda_k(1) &\iso& D^\lambda & 
			S^{\lambda^\top}_k(0)/S^{\lambda^\top}_k(1) &\iso& D^{\lambda^\top} \\ \\
			S^\lambda_k(1)/S^\lambda_k(2) &\iso& \bigoplus_{\mu\in r_\lambda \setminus \{ \lambda,\lambda^{\top M} \}} D^\mu  &
			S^{\lambda^\top}_k(1)/S^{\lambda^\top}_k(2) &\iso& \bigoplus_{\mu\in r_\lambda \setminus \{ \lambda, \lambda^{\top M} \}} D^{\mu^M} \\ \\
			S^\lambda_k(2)/S^\lambda_k(3) &\iso& D^{\lambda^{\top M}} & 
			S^{\lambda^\top}_k(2)/S^{\lambda^\top}_k(3) &\iso& D^{\lambda^{M}}
			\end{array}
		$$
		and all further layers are zero. Furthermore, $r_\lambda \setminus \{ \lambda,\lambda^{\top M} \} \neq \emptyset$, meaning
		all of the above layers are non-trivial.
	 \item[(ii)] If $\lambda$ is $p$-regular and $\lambda^\top$ is $p$-singular, then
		$$
			\begin{array}{rclrcl}
			S^\lambda_k(0)/S^\lambda_k(1) &\iso& D^\lambda  &
			S^{\lambda^\top}_k(0)/S^{\lambda^\top}_k(1) &\iso& 0 \\\\
			S^\lambda_k(1)/S^\lambda_k(2) &\iso& \bigoplus_{\mu\in r_\lambda \setminus \{ \lambda \}} D^\mu &
			S^{\lambda^\top}_k(1)/S^{\lambda^\top}_k(2) &\iso& \bigoplus_{\mu\in r_\lambda \setminus \{ \lambda \}} D^{\mu^M} \\\\
			S^\lambda_k(2)/S^\lambda_k(3) &\iso& 0 &
			S^{\lambda^\top}_k(2)/S^{\lambda^\top}_k(3) &\iso& D^{\lambda^{M}}
			\end{array}
		$$
		and all further layers are zero.
	 \item[(iii)] If $\lambda$ and $\lambda^\top$ are both $p$-singular, then there is a $p$-regular partition
		$\mu$ of $n$ (which will necessarily be the $p$-regularization of $\lambda$) such that
		$$
			\begin{array}{rclcrcl}
			S^\lambda_k(1)/S^\lambda_k(2) &\iso& D^\mu & \quad\quad &
			S^{\lambda^\top}_k(1)/S^{\lambda^\top}_k(2) &\iso& D^{\mu^M}
			\end{array}
		$$
		and all other layers are zero.
	\end{enumerate}
\begin{proof}
	By \cite[Theorem 4.8]{FayersJSFilt} (specialized to the defect two case) we have
	\begin{equation}\label{eqn_schaper_layer}
		S_k^\lambda(i)/S_k^\lambda(i+1) \iso \left( S_k^{\lambda^\top}(2-i)/S_k^{\lambda^\top}(3-i) \right) \otimes_\OO \sign
	\end{equation}
	In particular $S^\lambda_k(3)=\{0\}$.
	This clearly implies that the first and the third layer of the filtration are always as claimed. 
	Our claim on the middle layer in cases (i) and (ii) simply follows from the fact that all decomposition numbers are zero or one 
	(that is, any simple module that occurs as a composition factor of $S^\lambda_k$ does so with multiplicity one), and Remark \ref{remark_js_semisimple}.

	Now we show that when $\lambda$ and $\lambda^\top$ are both $p$-regular, the set $r_\lambda \setminus \{ \lambda,\lambda^{\top M} \}$ 
	is non-empty. Assume otherwise. By \cite[Corollary 13.18]{James} $S^\lambda_k$ is indecomposable, and thus 
	$\Ext^1_{k\Sym_n}(D^\lambda, D^{\lambda^{\top M}})$ must be non-zero. Now, as mentioned in Remark \ref{knownfacts}, the
	$\Ext$-quiver of a defect two block of a symmetric group is bipartite. The bipartition is given by the so-called 
	\emph{relative $p$-sign} (see \cite[Proposition 2.2.]{Def3Bipart}). Given any partition $\eta$, define its relative $p$-sign $\sigma_p(\eta)$ to be
	$(-1)^{\sum l_i}$, where $l_i$ are the leg lengths of a sequence of $p$-hooks that may be removed from
	$\eta$ to leave a $p$-core. Since for $p$ odd the leg length and the
	arm length of a $p$-hook always leave the same residue modulo two, we have $\sigma_p(\eta)=\sigma_p(\eta^\top)$ for 
	any partition $\lambda$. By \cite[Proposition 2.5.]{TanvDecomp}, for odd $p$ and $p$-regular $\eta$ of even weight,
	$\sigma_p(\eta)=\sigma_p(\eta^M)$ will hold. Therefore, $\sigma_p(\lambda^{\top M})=\sigma_p(\lambda)$, that is, 
	$D^\lambda$ and $D^{\lambda^{\top M}}$ are in the same part of the bipartition. But then, $\Ext^1$ between the two cannot
	be non-zero, giving us the desired contradiction.

	The only 
	part of our claim left to prove is that whenever $\lambda$ and $\lambda^\top$ are both $p$-singular,
	$S^\lambda_k$ will be simple. 
	By (\ref{eqn_schaper_layer}) it is clear that in this case, $S^\lambda_k \iso S_k^\lambda(1)/S_k^\lambda(2)$.
	According to Remark \ref{remark_js_semisimple} the module $S^\lambda_k$ is hence semisimple. But by 
	\cite[Corollary 13.18]{James}, $S^\lambda_k$ is also indecomposable. It follows that $S^\lambda_k$ is simple,
	as claimed. 
\end{proof}
\end{theorem}

\begin{remark}
	Note that the last remark  and theorem determine the submodule structure 
	of $S^\lambda_k$ for each partition $\lambda$ in a defect two block.
\end{remark}

\begin{lemma}\label{lemma_s_sdual}
	Let $\lambda$ be a partition of some $n$, and let $S^\lambda_K$ be equipped with the natural bilinear form
	inherited from $M^\lambda_K$. If $L\subset S^\lambda_K$ is a $\OO\Sym_n$-lattice, we denote its dual
	with respect to this form by $L^\sharp$. 
	Define $\hat{S}(j) := (p^{-j}\cdot S^\lambda_\OO) \cap S^{\lambda\sharp}_\OO$.
	Then there is an ascending filtration
	\begin{equation}
		S^\lambda_\OO = \hat{S}(0) \leq \hat{S}(1) \leq \ldots \leq \hat{S}(l) = S^{\lambda\sharp}_\OO \quad 
		\textrm{ for some $l \in \N$}
	\end{equation}
	and the quotients $\hat{S}(j)/\hat{S}(j-1)$ are isomorphic to $S^\lambda_k(j)$.
	\begin{proof}
		$$
			\hat{S}(j) / \hat{S}(j-1) \iso  \frac{p^{-j} S^\lambda_\OO \cap S^{\lambda\sharp}_\OO + p^{-j+1}S^\lambda_\OO}{p^{-j+1}S^\lambda_\OO} \iso \frac{S^\lambda_\OO\cap p^j S^{\lambda\sharp}_\OO + pS^\lambda_\OO}{pS^\lambda_\OO} = S^\lambda_k(j) 
		$$
	\end{proof}
\end{lemma}

\begin{theorem}\label{thm_grad_hull_def2}
	Let $\lambda$ be a $p$-regular partition in a defect two block. Let 
	$J_0, J_1$ and $J_2$ be the sets of $p$-regular partitions $\mu$ such that
	$D^\mu$ occurs in $S^\lambda_k(0)/S^\lambda_k(1)$, $S^\lambda_k(1)/S^\lambda_k(2)$ and $S^\lambda_k(2)/S^\lambda_k(3)$
	respectively. By $\eps^\lambda$ denote the primitive idempotent in $Z(K\Sym_n)$ belonging to $\lambda$. 
	Then the $\OO$-order $\eps^\lambda \OO\Sym_n$ is Morita-equivalent to the graduated order 
	$\Lambda = \Lambda(\OO,m,(1,\ldots,1))$ for an exponent matrix $m\in\Z_{\geq 0}^{r_\lambda \times r_\lambda}$ subject to the conditions:
	\begin{eqnarray}
	\label{exp_mat_def2_1}	m_{\alpha \lambda} &=& 0 \quad \forall \alpha \in J_0\cup J_1 \cup J_2 \\
	\label{exp_mat_def2_2}	m_{\lambda \alpha} &=& i \quad \forall \alpha\in J_i \quad \textrm{ for $i\in\{0,1,2\}$} \\
	\label{exp_mat_def2_3}	m_{\alpha\beta}-m_{\beta\alpha} &=& m_{\lambda\beta}-m_{\lambda\alpha}  \quad \forall \alpha,\beta \\
	\label{exp_mat_def2_4}	0 < m_{\alpha\beta}+m_{\beta\alpha} &\leq& 2 \quad\forall\alpha\neq\beta
	\end{eqnarray}
	These conditions completely determine the matrix $m$.

	Denote for each $\gamma\in r_\lambda$ by $e^\gamma$ the diagonal matrix unit in $\Lambda$ belonging 
	to $\gamma$. Then we may choose a Morita-equivalence $\mathcal{F}$ between $\modC_{\eps^\lambda \OO\Sym_n}$ and $\modC_\Lambda$ such that $\mathcal{F}(D^\gamma) \iso e^\gamma\Lambda / \Rad e^\gamma\Lambda$.
\begin{proof}
	Since $K$  and $k$ split $\Sym_n$ and all decomposition numbers are known to be $0$ or $1$ it follows that
	$\eps^\lambda \OO\Sym_n$ (as well as, of course, its basic algebra) is a graduated order (see Theorem \ref{thm_decomp_01}). Thus $\Lambda = \Lambda(\OO,m,(1,\ldots,1))$ for some exponent matrix $m$.
	We may assume without loss that
	$\mathcal{F}(S^{\lambda\sharp}_\OO)\iso \OO^{1\times J_0\cup J_1 \cup J_2}$. We may also assume
	that $\mathcal{F}(D^\gamma) \iso e^\gamma\Lambda / \Rad e^\gamma\Lambda$. At this point we have fixed an exponent matrix $m$, and we need to show that it satisfies (\ref{exp_mat_def2_1})-(\ref{exp_mat_def2_4}).
		
	The order $\eps^\lambda \OO\Sym_n$ carries an involution, as $K\Sym_n$ carries the standard 
	involution $g \mapsto g^{-1}$, and this involution fixes $\eps^\lambda$ and maps $\OO\Sym_n$ to itself.
	Dualizing followed by standard involution also fixes all simple modules (this fact is usually stated as ``The simple 
	$k\Sym_n$-modules are self-dual'').
	Hence the order $\Lambda$ may also be equipped with an involution ${^\circ}: \Lambda \rightarrow \Lambda$ that
	fixes all the $e^\gamma$ (this is by Corollary \ref{corollary_involution_basic_algebra} and Remark \ref{remark_involution}). By Theorem \ref{thm_sconst_invol} this implies (\ref{exp_mat_def2_3}).
	The equation (\ref{exp_mat_def2_4}) is just Remark \ref{remark_exp_mat_selfdual}.

	Since $S^\lambda_\OO$ has simple top $D^\lambda$, so does $\mathcal{F}(S^\lambda_{\OO})$. The uniqueness
	part of Remark \ref{remark_simple_graduated} thus implies $\mathcal{F}(S^\lambda_{\OO}) \iso e^\lambda\Lambda$. But $e^\lambda\Lambda \iso \left[ (p)_\OO^{m_{\lambda\alpha}} \right]_\alpha$, and therefore $m_{\lambda\alpha}$ equals (for each $\alpha$) the
	multiplicity of $\mathcal{F}(D^\alpha)$ in $\mathcal{F}(S_\OO^{\lambda\sharp})/\mathcal{F}(S_\OO^{\lambda})$, which has been determined in Lemma \ref{lemma_s_sdual}. This implies (\ref{exp_mat_def2_2}). Note that in principle $\mathcal{F}$ applied to an irreducible  lattice
	is, as a lattice in $K^{1\times r_\lambda}$, only determined up to multiplication by powers of $p$. 
	For the above quotient we choose however the maximal representative of $\mathcal{F}(S^\lambda_\OO)$ that is contained in $\OO^{1\times r_\lambda}
	= \mathcal{F}(S_\OO^{\lambda\sharp})$.
	
	We may equip the vector space $K^{1\times r_\lambda}$ with a $\Lambda$-equivariant non-degenerate 
	symmetric bilinear form (which one we choose is irrelevant for our purposes).
	Then for each $\OO\Sym_n$-lattice $L \leq S^\lambda_K$ we have $\mathcal{F}(L^\sharp)\iso \mathcal{F}(L)^\sharp$.
	This is best seen by choosing an involution-invariant idempotent in $e\in \eps^\lambda\OO\Sym_n$ that affords
	the Morita-equivalence, since 
	\begin{equation}
		L^\sharp \cdot e \iso \Hom_\OO(L,\OO)^\circ \cdot e \iso \Hom_\OO(L\cdot e^\circ, \OO)^\circ =\Hom_\OO(L\cdot e, \OO)^\circ \iso (L\cdot e)^\sharp
	\end{equation}

	By looking the standard bilinear pairing of $K^{1\times r_\lambda}$ and 
	$K^{r_\lambda \times 1}$ we see that 
	\begin{equation}
		\Hom_\OO(\mathcal{F}(S^{\lambda\sharp}_\OO),\OO) \iso_\Lambda \OO^{r_\lambda\times 1}
	\end{equation}
	On the other hand, as was just seen, 
	\begin{equation}
	\begin{array}{rcl}
	\Hom_\OO(\mathcal{F}(S^{\lambda\sharp}_\OO),\OO)
	&\iso& \Hom_\OO(\Hom_\OO(\mathcal{F}(S^{\lambda}_\OO), \OO)^\circ,\OO)\\& \iso& \mathcal{F}(S^\lambda_\OO)^\circ \iso (e^\lambda\Lambda)^\circ \iso \Lambda e^\lambda
	\end{array}
	\end{equation}
	This clearly implies (\ref{exp_mat_def2_1}).

	
	The conditions (\ref{exp_mat_def2_1})-(\ref{exp_mat_def2_4}) determine $m$, since 
	(\ref{exp_mat_def2_3}) determines for all $\alpha,\beta$ the difference  $m_{\alpha\beta}-m_{\beta\alpha}$, and
	hence determines $m_{\alpha\beta}+m_{\beta\alpha}$ modulo $2$. As (\ref{exp_mat_def2_4}) states
	that $m_{\alpha\beta}+m_{\beta\alpha}\in\{1,2\}$, this is already enough to determine the sum 
	$m_{\alpha\beta}+m_{\beta\alpha}$. This clearly determines the values of the $m_{\alpha\beta}$.
\end{proof}
\end{theorem}

\begin{remark}\label{remark_exp_mat_nonreglar}
	The preceding theorem determines the exponent matrices for every $\eps^\lambda\OO\Sym_n$ with
	$\lambda$ in a defect two block, even when $\lambda$ is $p$-singular. Namely, if $\lambda$ is $p$-singular, we have the following two cases:
	\begin{enumerate}
	\item \emph{$\lambda$ is $p$-singular and $\lambda^\top$ is $p$-regular:}\\
	We have an isomorphism 
	\begin{equation}
		\varphi:\ \eps^\lambda\OO\Sym_n \rightarrow \eps^{\lambda^\top}\OO\Sym_n: \ \eps^\lambda \cdot g 
		\mapsto \sign(g)\cdot \eps^{\lambda^\top}\cdot g
	\end{equation}
	and when we retract the simple $\eps^{\lambda^\top}\OO\Sym_n$-modules with $\varphi$, we get
	$\varphi_*(D^{\mu^M}) \iso D^\mu$. Hence 
	\begin{equation}
		m^\lambda_{\mu\nu} = m^{\lambda^\top}_{\mu^M\nu^M} \quad \textrm{ for all } \mu,\nu \in r_\lambda
	\end{equation}
	where $m^\lambda$ and $m^{\lambda^\top}$ denote the exponent matrices of $\eps^\lambda\OO\Sym_n$ and
	$\eps^{\lambda^\top}\OO\Sym_n$. $m^{\lambda^\top}$ has of course been determined by the last theorem.
	\item \emph{$\lambda$ is $p$-singular and $\lambda^\top$ is $p$-singular:}\\
	According to Theorem \ref{thm_structure_schaper} the set $r_\lambda$ will contain just one element,
	and hence the exponent matrix will be the $1\times 1$ zero matrix.
	\end{enumerate}
\end{remark}

\begin{corollary}\label{corollary_ext_def2}
Let $\lambda$ be a partition in a defect two block. Then the $\Ext$-quiver of
$k\otimes_\OO \eps^\lambda\OO\Sym_n$ is maximally bipartite. More precisely, this means that there is an edge between 
any vertex pertaining to a constituent of $S^\lambda_k(1)/S^\lambda_k(2)$ and any vertex pertaining to a 
constituent of either $S^\lambda_k(0)/S^\lambda_k(1)$ or  $S^\lambda_k(2)/S^\lambda_k(3)$.
\begin{proof}
	By Remark \ref{remark_exp_mat_nonreglar} we may assume that $\lambda$ is $p$-regular (and we do so in what follows). We adopt the notation 
	of Theorem \ref{thm_grad_hull_def2}. 
	The group $\Sym_{J_1}$ acts naturally via automorphisms 
	on the basic order of $\eps^\lambda \OO\Sym_n$. In particular it acts via quiver automorphisms
	on the $\Ext$-quiver of $k\otimes_\OO \eps^\lambda\OO\Sym_n$ by permuting the vertices labeled
	by elements of $J_1$. This can easily be derived from the fact that the set of equations 
	(\ref{exp_mat_def2_1})-(\ref{exp_mat_def2_4}) is invariant  under the operation of $\Sym_{J_1}$ on the indices, and those equations
	determine the exponent matrix $m$ completely.

	$J_0$ and $J_2$ each consist of at most one partition. First suppose that 
	$J_1 \neq \emptyset$. Then by Theorem \ref{thm_structure_schaper} the top of $S^\lambda_k$ has a single constituent
	labeled by the partition in $J_0$, and the socle of $S^\lambda_k$ has constituents labeled by the partitions
	in $J_2$ (also, of course, at most one). Therefore $S^\lambda_k/\Soc S^\lambda_k$ has at least one non-semisimple quotient of length two, 
	implying the existence of an edge from the partition in $J_0$ to one partition in $J_1$. Provided 
	$J_2\neq\emptyset$, the module $\Rad S^\lambda_k$ has a non-semisimple submodule of length two, implying the existence of an edge from the element of $J_2$ to one element of $J_1$. Now using the action of $\Sym_{J_1}$
	we conclude that the $\Ext$-quiver has at least the postulated  edges. The case $J_1=\emptyset$ is trivial, since then by Theorem \ref{thm_structure_schaper} the set $J_2$ is also empty, that is, the $\Ext$-quiver consists of only a singe vertex. 

	As we already mentioned in Remark \ref{knownfacts}, the $\Ext$-quiver of any defect two block
	of a symmetric group is known to be bipartite by \cite{ChuangTanYoung}. 
	We can use the epimorphism $k\Sym_n \twoheadrightarrow k\otimes_\OO\eps^\lambda\OO\Sym_n$ to retract modules and sequences of modules, in particular
	simple modules and extensions of simple modules. Hence the $\Ext$-quiver of $k\otimes_\OO\eps^\lambda\OO\Sym_n$ is a
	sub-quiver of the bipartite  $\Ext$-quiver of the defect two block. It will therefore be bipartite as well. But if any further edges were to be
	added to the quiver constructed above, it would cease to be bipartite (for then there would be a closed path
	of length three). Hence we have constructed the full $\Ext$-quiver of $k\otimes_\OO\eps^\lambda \OO\Sym_n$.
\end{proof}
\end{corollary}

\begin{lemma}\label{lemma_bipartit_schaper}
	Let $\lambda$ be a 
	partition in a defect two block, and let $i \in\Z_{\geq 0}$. 
	If $D^\gamma$ is a simple module occurring in $S^\lambda_k(i)/S^\lambda_k(i+1)$ and 
	$D^\omega$ is a simple module occurring in $S^\lambda_k(i+1)/S^\lambda_k(i+2)$, then
	$\Ext^1_{k\Sym_n}(D^\gamma, D^\omega)\neq \{0\}$. 
	On the other hand, whenever two simple modules $D^\gamma$ and $D^\omega$ occur in the same layer of the Jantzen-Schaper-filtration,
	then $\Ext^1_{k\Sym_n}(D^\gamma, D^\omega)=\{0\}$.
\begin{proof}
	This follows directly from Corollary \ref{corollary_ext_def2}.
\end{proof}
\end{lemma}

%

\begin{theorem}\label{thm_ext_cartan}
	Let $\lambda$ and $\mu$ be two distinct $p$-regular partitions in some defect two block.
	If $\Ext_{k\Sym_n}^1(D^\lambda, D^\mu) \neq \{0\}$, then both of the following hold:
	\begin{enumerate}
	 \item[(i)]  $|c_\lambda \cap c_\mu| = 2$
	 \item[(ii)] $\lambda \in c_\lambda\cap c_\mu$ or $\mu \in c_\lambda\cap c_\mu$
	\end{enumerate}
\begin{proof}
	To prove the claim of (i), we argue by contradiction. We know that $|c_\lambda \cap c_\mu|\leq 2$ by Remark \ref{knownfacts} (i). So let $c_\lambda \cap c_\mu$ consist of just
	one element, say $\eta$. Then by Theorem \ref{thm_structure_schaper} and Lemma \ref{lemma_bipartit_schaper}, 
	$D^\lambda$ and $D^\mu$ occur in successive Jantzen-Schaper layers of $S^\eta_k$. Hence, by Theorem
	\ref{thm_grad_hull_def2}, $m_{\mu\lambda} + m_{\lambda\mu} = 1$ (where $m$ is the exponent matrix
	of $\eps^\eta \OO\Sym_n$).
 	 But if $e_\lambda$ and $e_\mu$ are primitive idempotents in $\OO\Sym_n$
 	corresponding to $D^\lambda$ and $D^\mu$, then 
 	$e_\lambda \OO\Sym_n e_\mu \OO\Sym_n e_\lambda = \left\langle p \cdot \eps^\eta \cdot e_\lambda \right\rangle_\OO$
	(this is easy to see if one identifies $\eps^\eta \OO\Sym_n$ with $\Lambda(\OO,m,d)$ for the appropriate dimension vector $d$, and assumes without loss that $\eps^\eta e_\lambda$ and $\eps^\eta e_\mu$
	are equal to diagonal matrix units in $\Lambda(\OO,m,d)$). $\OO\Sym_n$
	is a self-dual (and so in particular integral) lattice  with respect to the bilinear form 
	$T:\ (a,b)\mapsto \frac{1}{n!}\sum_{\varphi} \chi^\varphi(1)\chi^\varphi(ab)$ (where $\chi^\varphi$ is the irreducible character associated to the Specht module $S^\varphi_{K}$).
	Now $T(p \cdot \eps^\eta \cdot e_\lambda,1)=p\cdot T(\eps^\eta \cdot e_\lambda,1)=p\cdot \frac{\chi^\eta(1)}{n!}\cdot\chi^\eta(e_\lambda)$. However $\nu_p(p)+\nu_p(\frac{\chi^\eta(1)}{n!})+\nu_p(\chi^\eta(e_\lambda)) = 1 - 2 + 0 = -1$, and thus 
	$T(p \cdot \eps^\eta \cdot e_\lambda,1)\notin \OO$, in contradiction to the integrality of $\OO\Sym_n$. 

	Now we prove (ii).
	Let $\nu$ be an element of $c_\lambda \cap c_\mu$.
	By Lemma \ref{lemma_bipartit_schaper}, either $\lambda$ or $\mu$ must occur in one of $S^\nu_k(0)/S^\nu_k(1)$
	or $S^\nu_k(2)/S^\nu_k(3)$. By Theorem \ref{thm_structure_schaper} it follows that
	$\nu \in \{ \lambda, \mu, \lambda^{M\top}, \mu^{M\top}\}$. Since we know already that $|c_\lambda\cap c_\mu|=2$, we only have to check that
	$c_\lambda \cap c_\mu \neq \{ \lambda^{M\top}, \mu^{M\top} \}$.
	Suppose the contrary. Then $S_k^{\lambda^{M\top}}$ has $D^\mu$ as a composition factor, and 
	therefore $S_k^{\lambda^M}$ has $D^{\mu^M}$ as a composition factor. It follows $\mu^M\triangleright\lambda^M$.  But in the same way the fact that $S_k^{\mu^{M\top}}$ has $D^\lambda$ as a composition factor implies that $\lambda^M\triangleright\mu^M$, which yields the desired contradiction.
\end{proof}
\end{theorem}

At this point we fix a defect two block $B$ of some $\OO\Sym_n$, and we wish to describe its basic order, which we 
shall denote by $\Lambda$. We describe $\Lambda$ as an order in the $K$-algebra $A$ which we are about to define.

\begin{defi}
	The exponent matrices for $B$ determined in Theorem \ref{thm_grad_hull_def2} and Remark \ref{remark_exp_mat_nonreglar} shall be denoted by
	$m^\lambda_{\mu\nu}$. By $d_\lambda$ we denote the dimension of the Specht module $S^\lambda_K$.
\end{defi}

\begin{defi}
	Define a $K$-algebra $A$ spanned by a $K$-basis
	\begin{equation}
		\eps_{\mu\nu}^{\lambda} \quad \textrm{ for $\lambda$ a partition in $B$ and $\mu,\nu\in r_\lambda$}
	\end{equation}
	equipped with the following multiplication law:
	\begin{equation}
		\eps_{\mu\nu}^{\lambda} \cdot \eps_{\tilde{\mu}\tilde{\nu}}^{\tilde{\lambda}}
		= \delta_{\lambda\tilde{\lambda}} \cdot \delta_{\nu\tilde{\mu}} \cdot \eps_{\mu\tilde{\nu}}^{\lambda}
	\end{equation}
	Note that this $A$ is isomorphic to a direct sum of full matrix algebras over $K$, and the $\eps_{\mu\nu}^{\lambda}$ are just the matrix units. Note also that $A \iso K\otimes_\OO\Lambda$.
	We will henceforth assume that $\Lambda$ is embedded in $A$.
	
	The central primitive idempotents in $A$ are given by
	\begin{equation}
		\eps^\lambda  := \sum_{\mu \in r_\lambda} \eps_{\mu\mu}^{\lambda}
	\end{equation}
	and we shall assume without loss that for each $\lambda$
	\begin{equation}\label{ass_proj_surj_def2}
		\eps^\lambda\Lambda = \bigoplus_{\mu,\nu\in r_\lambda} 
		\langle p^{m^\lambda_{\mu\nu}}\cdot \eps_{\mu\nu}^{\lambda} \rangle_\OO
	\end{equation}
	and that for each $p$-regular $\mu$ the idempotent $\sum_{\lambda\in c_\mu} \eps_{\mu\mu}^{\lambda}$
	is a primitive idempotent in $\Lambda$ (corresponding to $D^\mu$). This can all be achieved by conjugation within $A$.
\end{defi}

\begin{theorem}\label{thm_basic_order}
	The order $\Lambda$ is conjugate in $A$ to the $\OO$-algebra generated by the following elements of $A$:
	For each $p$-regular $\mu$ in $B$ the idempotent
	\begin{equation}
		e_\mu := \sum_{\lambda \in c_\mu} \eps_{\mu\mu}^{\lambda}
	\end{equation}
	and for each (ordered) pair $(\mu,\nu)$ of (distinct) $p$-regular partitions in $B$ with $\Ext^1_{k\otimes_\OO B}(D^\mu,D^\nu)\neq\{0\}$ an element $x_{\mu\nu}$, which is defined as 
	\begin{equation}
		x_{\mu\nu} := p^{m^\lambda_{\mu\nu}}\cdot \eps_{\mu\nu}^{\lambda} + p^{m^\eta_{\mu\nu}}\cdot \eps_{\mu\nu}^{\eta}
		\quad\textrm{where } c_\mu\cap c_\nu = \{\lambda, \eta\} 
	\end{equation}
	if $ \mu > \nu$, respectively
	\begin{equation}
	 	x_{\mu\nu} := p^{m^\lambda_{\mu\nu}}\cdot \eps_{\mu\nu}^{\lambda} - \frac{d_\lambda}{d_\eta}\cdot p^{m^\eta_{\mu\nu}}\cdot \eps_{\mu\nu}^{\eta}
		\quad\textrm{ where } c_\mu\cap c_\nu = \{\lambda, \eta\}
	\end{equation}
	if $\mu < \nu$.
\begin{proof}
	It follows easily from Nakayama's lemma (for $\OO$-modules) that a set of elements of $\Lambda$ 
	generate $\Lambda$ as an $\OO$-algebra if and only if their images in $k\otimes_\OO\Lambda$ generate $k\otimes_\OO\Lambda$
	as a $k$-algebra.
	It is well known (see for instance \cite[Proposition 4.1.7]{Benson}, and be aware that the condition
	``$k$ algebraically closed'' may be replaced by ``$k$ is a splitting field'') that 
	a full set of primitive idempotents $\overline{e}_\mu$ (with the natural choice of indices) and any basis of $\overline{e}_\mu \Jac(k\otimes_\OO\Lambda)/\Jac(k\otimes_\OO\Lambda)^2 \overline{e}_\nu$
	(where $\mu,\nu$ run over all $p$-regular partitions) generates $k\otimes_\OO \Lambda$. 
	By \cite[Proposition 2.4.3]{Benson} we have
	\begin{equation}
	 \dim_k \overline{e}_\mu\cdot \left(\Jac(k\otimes_\OO\Lambda)/\Jac(k\otimes_\OO\Lambda)^2\right)\cdot \overline{e}_\nu = \dim_k \Ext^1_{k\otimes_\OO B}(D^\mu,D^\nu)
	\end{equation}
	Moreover we know (as mentioned in Remark \ref{knownfacts}) that all $\Ext^1_{k\otimes_\OO B}(D^\mu,D^\nu)$ are at most one-dimensional.

	Now pick a specific pair $\mu,\nu$ of $p$-regular partitions with $\mu > \nu$
	such that $\Ext^1_B(D^\mu,D^\nu)\neq\{0\}$. Our goal is to pick some element in
	$e_\mu\Lambda e_\nu$ that is suitable as a generator due to the above considerations.

	First we should note that since $\Lambda$ is a self-dual order with respect to the
	symmetric $\Lambda$-equivariant bilinear form
	\begin{equation}
		A \times A \rightarrow K: (a,b) \mapsto \frac{1}{n!}\sum_{\lambda}d_\lambda\cdot \Tr (\eps^\lambda \cdot a \cdot b)
	\end{equation}
	we can define the bilinear pairing
	\begin{equation}\label{bifo_lambdamunu}
		\begin{array}{rccc}
		T: & e_\nu A e_\mu \times e_\mu A e_\nu &\rightarrow& K\\
		&  \left(\sum_{\lambda \in c_\mu\cap c_\nu} 
		f_\lambda \cdot \eps_{\nu\mu}^{\lambda},\ \sum_{\lambda \in c_\mu\cap c_\nu} 
		g_\lambda \cdot \eps_{\mu\nu}^{\lambda} \right) &\mapsto& \frac{1}{n!}\cdot\sum_{\lambda\in c_\mu\cap c_\nu} d_\lambda\cdot f_\lambda\cdot g_\lambda 
		\end{array}
	\end{equation}
	to get $e_\nu \Lambda e_\mu = \left\{ v \in e_\nu A e_\mu\ |\ T(v,e_\mu \Lambda e_\nu) \subseteq \OO \right\}$ (and the analogous equation for $e_\mu \Lambda e_\nu$).
	We will use this together with the fact that $\nu_p\left(\frac{d_\lambda}{n!}\right) = -2$ for all
	$\lambda$.
	We distinguish the following cases:
	\begin{enumerate}
	 \item[(i)] $c_\mu\cap c_\nu = \{\lambda,\eta\}$, $m^\lambda_{\mu\nu}+m^\lambda_{\nu\mu}=2$ and
		$m^\eta_{\mu\nu}+m^\eta_{\nu\mu}=2$. By Theorem \ref{thm_structure_schaper} and
		Theorem \ref{thm_grad_hull_def2} this could only happen if $\{\lambda,\eta\} = \{\mu,\nu\}$. 
		But then $D^\mu$ is a composition factor of $S^\nu$, so $\mu\triangleright \nu$, and $D^\nu$ is a composition factor of 
		$S^\mu$, so $\nu\triangleright\mu$. Clearly this is a contradiction, so this case does not occur at all.
	 \item[(ii)] $c_\mu\cap c_\nu = \{\lambda,\eta\}$ and $m^\lambda_{\mu\nu}+m^\lambda_{\nu\mu}=1$. 
		In this case $T(e_\nu \Lambda e_\mu, p^{m^\lambda_{\mu\nu}}\cdot \eps_{\mu\nu}^{\lambda})= p^{-1}\cdot\OO$, which implies $p^{m^\lambda_{\mu\nu}}\cdot\eps_{\mu\nu}^{\lambda}\notin e_\mu\Lambda e_\nu$
		(note however that $p^{m^\lambda_{\mu\nu}+1}\cdot\eps_{\mu\nu}^{\lambda}$ is in $e_\mu \Lambda e_\nu$ by the same argument), and thus
		\begin{equation}e_\mu\Lambda e_\nu \subsetneqq \langle p^{m^\lambda_{\mu\nu}} \cdot \eps_{\mu\nu}^{\lambda}\rangle_\OO \oplus \langle p^{m^\eta_{\mu\nu}} \cdot \eps_{\mu\nu}^{\eta}\rangle_\OO\end{equation} However the projection onto 
		each summand (that is, multiplication by $\eps^\lambda$) has to be surjective by (\ref{ass_proj_surj_def2}), and so there is an element in $e_\mu \Lambda e_\nu$ of the form 
		$p^{m^\lambda_{\mu\nu}}\cdot \eps_{\mu\nu}^{\lambda}+\alpha_{\mu\nu}^{\eta} \cdot p^{m^\eta_{\mu\nu}}\cdot \eps_{\mu\nu}^{\eta}$ for some $\alpha_{\mu\nu}^{\eta} \in \OO^{\times}$.
		So we can state that
		\begin{equation}
			e_\mu \Lambda e_\nu = \left\langle  p^{m^\lambda_{\mu\nu}}\cdot \eps_{\mu\nu}^{\lambda}
			+  \alpha_{\mu\nu}^{\eta}\cdot p^{m^\eta_{\mu\nu}}\cdot \eps_{\mu\nu}^{\eta},\ 
			 p^{m^\lambda_{\mu\nu}+1}\cdot \eps_{\mu\nu}^{\lambda}  \right\rangle_\OO
		\end{equation}
		and by dualizing it follows that
		\begin{equation}
			e_\nu \Lambda e_\mu = \left\langle  p^{m^\lambda_{\nu\mu}}\cdot \eps_{\nu\mu}^{\lambda}
			-  \frac{d_\lambda}{d_\eta}\cdot(\alpha_{\mu\nu}^{\eta})^{-1}\cdot p^{m^\eta_{\nu\mu}}\cdot \eps_{\nu\mu}^{\eta},\ 
			 p^{m^\lambda_{\nu\mu}+1}\cdot \eps_{\nu\mu}^{\lambda} \right\rangle_\OO
		\end{equation}
		 Theorem \ref{thm_amal_comm} and Remark \ref{remark_amal_comm} imply that
		\begin{equation}
			(\eps^\lambda+\eps^\eta)\cdot e_\mu \Lambda e_\mu = \langle 
				\eps_{\nu\mu}^{\lambda} + \eps_{\nu\mu}^{\eta}, p \cdot \eps_{\nu\mu}^{\lambda}
			 \rangle_\OO
		\end{equation}
		and therefore $\langle p^{m^\lambda_{\nu\mu}+1}\cdot \eps_{\nu\mu}^{\lambda},
		p^{m^\eta_{\nu\mu}+1}\cdot \eps_{\nu\mu}^{\eta}  \rangle_\OO$ is the unique maximal
		$e_\mu \Lambda e_\mu$-submodule of $e_\nu \Lambda e_\mu = e_\nu \Jac(\Lambda) e_\mu$. It is therefore equal to
		$e_\mu \Jac(\Lambda)^2 e_\nu$. Hence we may take 
		\begin{equation}\label{eqn_jklj38847j}x_{\mu\nu} =  p^{m^\lambda_{\mu\nu}}\cdot \eps_{\mu\nu}^{\lambda}
			+  \alpha_{\mu\nu}^{\eta}\cdot p^{m^\eta_{\mu\nu}}\cdot \eps_{\mu\nu}^{\eta}
		\end{equation}
		as a generator (since it is not contained in $e_\mu \Jac(\Lambda)^2 e_\nu$), and by the same argument we may pick
		\begin{equation}\label{eqn_jokj748764bgv}
			x_{\nu\mu}=p^{m^\lambda_{\nu\mu}}\cdot \eps_{\nu\mu}^{\lambda}
			- \frac{d_\lambda}{d_\eta}\cdot(\alpha_{\mu\nu}^{\eta})^{-1}\cdot p^{m^\eta_{\nu\mu}}\cdot \eps_{\nu\mu}^{\eta}		 		
		\end{equation}
	\end{enumerate}
	Now we have to show that all the $\alpha_{\mu\nu}^{\eta}$ may be chosen to be equal to one.
	To do that, first note that due to Theorem \ref{thm_ext_cartan} we may assume that all parameters
	are of the form $\alpha_{\mu\nu}^{\nu}$ (for $p$-regular partitions $\mu > \nu$). Of course, in this case, it may 
	also be assumed that $\nu$ is the lexicographically greatest element in $c_\mu\cap c_\nu$.

	Assume $\nu_0$ is a $p$-regular partition in $B$ such that all $\alpha_{\mu\nu}^{\nu}$ with $\nu < \nu_0$ are equal to one and $\alpha_{\mu\nu_0}^{\nu_0}\neq 1$ for some $\mu$.
	Assume moreover that $\nu_0$ is lexicographically maximal with respect to this property.
	Our goal is to show that after conjugation by an appropriate unit in $A$ and renormalization of the
	generators by multiplying with elements of $\OO^\times$ afterwards, we can make it so that all $\alpha_{\mu\nu}^{\nu}$ with $\nu \leq \nu_0$ equal one, which yields that without loss, 
	all $\alpha^\nu_{\mu\nu}$ may be chosen equal to one.
	To do this, we conjugate with a $u$ (i. e., replace each $x_{\mu\nu}$ by $u^{-1}\cdot x_{\mu\nu}\cdot u$), where 
	\begin{equation}
		u := \sum_{\mu \in r_{\nu_0}} \alpha_{\mu\nu_0}^{\nu_0} \cdot \eps_{\mu\mu}^{\nu_0} + \sum_{\lambda \neq \nu_0} \eps^\lambda \quad \in A^\times
	\end{equation}
	Note that in this formula we take those $\alpha_{\mu\nu_0}^{\nu_0}$ that are not defined to equal one. The conjugation with this unit
	will obviously make all $\alpha_{\mu\nu_0}^{\nu_0}$ equal to one, and not affect any $\alpha_{\mu\nu}^{\nu}$
	with $\nu < \nu_0$ (since all elements of $c_\mu\cap c_\nu$ will be lexicographically smaller than $\nu_0$).
	After renormalizing the other generators that were altered by the conjugation (to make them look as in (\ref{eqn_jklj38847j}) respectively (\ref{eqn_jokj748764bgv}) again), we have $\alpha^\nu_{\mu\nu}=1$ for all $\nu\leq\nu_0$. That concludes the proof. 
\end{proof}
\end{theorem}

It has been proved in \cite[Corollary 5.4.5.]{Peach} that defect two blocks of symmetric groups over $k$ are tightly graded. The following reproves that result (in a very simple fashion), and additionally shows that the 
images of the $x_{\mu\nu}$ (as defined in the theorem above) in the basic algebra over $k$ are homogeneous generators. That should in particular simplify the calculation of the quiver relations from our description of the block. 

\begin{corollary}
	Let $\Lambda = \OO\langle \{e_\mu\}_\mu, \{x_{\mu\nu}\}_{\mu,\nu}\rangle$ as in Theorem 
	\ref{thm_basic_order}. Let $Q$ be the $\Ext$-quiver, denote by $E_\mu$ the vertices
	and denote by $X_{\mu\nu}$ an edge from $E_\mu$ to $E_\nu$. By $kQ$ we denote the quiver algebra
	(with multiplication convention $X_{\mu\nu}\cdot X_{\nu\tau} \neq 0$). Then the kernel of the epimorphism
	\begin{equation}
		\Phi:\ kQ \twoheadrightarrow \Lambda / p\Lambda: \ \left\{ \begin{array}{c}X_{\mu\nu} \mapsto x_{\mu\nu} + p\Lambda\\ E_\mu \mapsto e_\mu + p\Lambda \end{array}\right.
	\end{equation}
	is a homogeneous ideal, where we define the vertices of $Q$ to be homogeneous of degree zero and the arrows to
	be homogeneous of degree one.
	\begin{proof}
		Since $\Ker \Phi = \bigoplus_{\mu,\nu} E_\mu \cdot \Ker \Phi \cdot  E_\nu$, and a path connecting
		$E_\mu$ with $E_\nu$ has even respectively odd length if and only if $\mu$ and $\nu$ lie in the same part respectively
		in different parts of the bipartition, we may assume that $\Ker \Phi$ is generated by elements 
		that involve only paths of even length and elements that only involve paths of odd length.
		By general theory we may assume that all paths involved in any element of $\Ker \Phi$ have at least
		length two. By \cite[Theorem I]{ScopesDef2}, the projective indecomposables of $\Lambda/p\Lambda$
		have common Loewy length five, i. e. $\Phi$ maps every path of length $\geq 5$  to zero. 
		A path of length four will correspond to a top onto socle endomorphism of a projective indecomposable.
		Thus all paths of length four that start and end at a separate vertex are sent to zero under 
		$\Phi$, and all paths of length four that start and end at the same fixed vertex of $Q$ will
		be mapped by $\Phi$ into a one-dimensional subspace of $\Lambda/p\Lambda$.

		So, considering all of this, all we need to show is that if $Y_\mu := X_{\mu\alpha}X_{\alpha\beta}X_{\beta\gamma}X_{\gamma\mu}$ is not in $\Ker \Phi$, then
		neither is $Y_\mu + \sum_\nu q_\nu \cdot X_{\mu\nu}X_{\nu\mu}$ for any choice of $q_\nu\in k$.
		It follows easily from Theorem \ref{thm_ext_cartan} that $|c_\mu \cap c_\alpha \cap c_\beta \cap c_{\gamma}|=1$, and let us denote the single element of this set by $\lambda$. Hence $y_\mu := x_{\mu\alpha}x_{\alpha\beta}x_{\beta\gamma}x_{\gamma\mu}$
		is equal to $v \cdot \eps^\lambda_{\mu\mu}$ for some $v \in \OO$. The fact that
		$y_\mu \in \Lambda \setminus p\Lambda$ implies $\nu_p(v)=2$. Let $T: e_\mu\Lambda e_\mu \times
		e_\mu\Lambda e_\mu \rightarrow \OO$ the symmetric bilinear form on $e_\mu \Lambda e_\mu$ as given 
		in  (\ref{bifo_lambdamunu}). Then $T(y_\mu,1)\in \OO^\times$. On the other hand
		$T(x_{\mu\nu}x_{\nu\mu},1)=0$ for any $\nu$ by definition of the $x_{\mu\nu}$. Hence
		$T(y_\mu + \sum_\nu \hat{q}_\nu \cdot x_{\mu\nu}x_{\nu\mu}, 1)\in \OO^\times$ for any choice 
		of $\hat{q}_\nu \in \OO$, which implies $y_\mu + \sum_\nu \hat{q}_\nu \cdot x_{\mu\nu}x_{\nu\mu}\notin p\Lambda$. Thus $\Phi(Y_\mu + \sum_\nu q_\nu \cdot X_{\mu\nu}X_{\nu\mu})\neq 0$.
	\end{proof}
\end{corollary}

\newcommand{\pone}{\lambda}
\newcommand{\ptwo}{\mu}
\newcommand{\pthree}{\nu}
\newcommand{\pfour}{\eta}
\newcommand{\pfive}{\varphi}
\newcommand{\psix}{{\tilde{\pfour}}}
\newcommand{\pseven}{\tilde{\pthree}}
\newcommand{\peight}{\tilde{\ptwo}}
\newcommand{\pnine}{\tilde{\pone}}

\begin{example}
	We look at the principal block of $\Z_3\Sym_7$. The decomposition matrix is given as follows
	(we assign arbitrary names to the partitions in order to unclutter notation a bit):
	\begin{equation}
		\begin{array}{ccrccccc}
			 \textrm{Dim.} & \textrm{Name} & & (7) & (5,2) & (4,3) & (4,2,1) & (3,2,1^{2}) \\ 
			\cmidrule[1.2pt]{1-8}
 			1 & \pone &(7) & 1 & . & . & . & .  \\ 
			14 & \ptwo & (5,2) & 1 & 1 & . & . & .  \\ 
			14 & \pthree & (4,3) & . & 1 & 1 & . & .  \\ 
			35 & \pfour & (4,2,1) & 1 & 1 & 1 & 1 & .  \\ 
			20 & \pfive & (4,1^{3}) & . & . & . & 1 & .  \\ 
			35 & \psix & (3,2,1^{2}) & 1 & . & 1 & 1 & 1  \\ 
			14 & \pseven & (2^{3},1) & 1 & . & . & . & 1  \\ 
			14 & \peight &(2^{2},1^{3}) & . & . & 1 & . & 1  \\ 
			1 & \pnine & (1^{7}) & . & . & 1 & . & .  
		\end{array}
	\end{equation}
	and the $\Ext$-quiver is given by
	$$
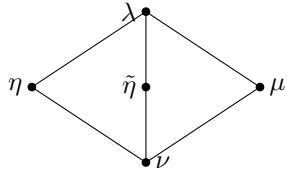
\begin{tikzpicture}[x=1cm,y=1cm]
		\draw (-0,-1) -- (-1.5,0);
		\draw (-0,1) -- (-1.5,0);
		\draw (-0,-1) -- (0,0);
		\draw (-0,1) -- (0,0);
		\draw (-0,-1) -- (1.5,0);
		\draw (-0,1) -- (1.5,0);
		\draw[fill=black] (-0,1) circle (0.05) node [left] (A)  {$\pone$};
		\draw[fill=black] (-1.5,0) circle (0.05) node[left] (B) {$\pfour$};
		\draw[fill=black] (0,0) circle (0.05) node[left] (C) {$\psix$};
		\draw[fill=black] (1.5,0) circle (0.05) node[right] (E) {$\ptwo$};
		\draw[fill=black] (-0,-1) circle (0.05) node[right] (F) {$\pthree$};
	\end{tikzpicture}$$
	The $4\times 4$-exponent matrices are given as follows
	\begin{equation}
		m^{\pfour}=m^{\psix}= \left( \begin{array}{cccc} 0&0&0&0\\1&0&1&0\\1&1&0&0\\2&1&1&0 \end{array}\right)
	\end{equation}
	with row/column indexing $\left(\ptwo,\pone,\pthree,\pfour\right)$ and $\left(\pfour,\pone,\pthree,\psix\right)$.
	The $2\times 2$-exponent matrices are
	\begin{equation}
		m^{\ptwo}=m^{\pthree}=m^{\pseven}=m^{\peight}=\left(\begin{array}{cc} 0 & 1 \\ 0 & 0 \end{array}\right)
	\end{equation}
	with rwo/column indexing $\left( \pone,\ptwo \right)$, $\left( \ptwo,\pthree \right)$, $\left( \pone,\psix \right)$
	and $\left( \pthree,\psix \right)$. The $1\times 1$-exponent matrices are of course trivial.

	By Theorem \ref{thm_basic_order} we now get the following generators for the basic order:
	\begin{equation}
		\scriptstyle{
		\begin{array}{rclrcl}
		x_{\pone \pfour}&=& \eps^{\pfour}_{\pone \pfour} + 3\cdot \eps^{\psix}_{\pone \pfour} 
		& 
		x_{\pfour \pone}&=& 3\cdot \eps^{\pfour}_{\pfour \pone} -\eps^{\psix}_{\pfour \pone}  \\ 
		x_{\pone \psix}&=& \eps^{\psix}_{\pone \psix} + 3\cdot \eps^{\pseven}_{\pone \psix} 
		& 
		x_{\psix \pone}&=&  3\cdot \eps^{\psix}_{\psix \pone} - \eps^{\pseven}_{\psix \pone} \\ 
		x_{\pone \ptwo}&=& 3\cdot \eps^{\ptwo}_{\pone \ptwo} + 3\cdot \eps^{\pfour}_{\pone \ptwo} 
		& 
		x_{\ptwo \pone}&=&   \eps^{\ptwo}_{\ptwo \pone} - \eps^{\pfour}_{\ptwo \pone} \\ 
		x_{\pthree \pfour}&=& \eps^{\pfour}_{\pthree \pfour} + 3\cdot \eps^{\psix}_{\pthree \pfour} 
		& 
		x_{\pfour \pthree}&=&3\cdot \eps^{\pfour}_{\pfour \pthree} -  \eps^{\psix}_{\pfour \pthree}  \\
		x_{\pthree \psix}&=& \eps^{\psix}_{\pthree \psix} + 3\cdot \eps^{\peight}_{\pthree \psix} 
		& 
		x_{\psix \pthree}&=&  3\cdot \eps^{\psix}_{\psix \pthree} - \eps^{\peight}_{\psix \pthree}  \\
		x_{\ptwo \pthree}&=& 3\cdot \eps^{\pthree}_{\ptwo \pthree} + \eps^{\pfour}_{\ptwo \pthree} 
		& 
		x_{\pthree \ptwo}&=& \eps^{\pthree}_{\pthree \ptwo} - 3\cdot \eps^{\pfour}_{\pthree \ptwo} 
		\\
		\end{array}}
	\end{equation}
	and of course the following idempotents:
	\begin{equation}
		\scriptstyle{
		\begin{array}{rcl}		
		e_{\pone} &=& \eps^{\pone}_{\pone \pone} + \eps^{\ptwo}_{\pone \pone} + \eps^{\pfour}_{\pone \pone}
			+ \eps^{\psix}_{\pone \pone} + \eps^{\pseven}_{\pone \pone} \\
		e_{\ptwo} &=& \eps^{\ptwo}_{\ptwo \ptwo} + \eps^{\pthree}_{\ptwo \ptwo} + \eps^{\pfour}_{\ptwo \ptwo} \\
		e_{\pthree} &=& \eps^{\pthree}_{\pthree \pthree} + \eps^{\pfour}_{\pthree \pthree} + \eps^{\psix}_{\pthree \pthree}
			+ \eps^{\peight}_{\pthree \pthree} + \eps^{\pnine}_{\pthree \pthree} \\
		e_{\pfour} &=& \eps^{\pfour}_{\pfour \pfour} + \eps^{\pfive}_{\pfour \pfour} + \eps^{\psix}_{\pfour \pfour} \\
		e_{\psix} &=& \eps^{\psix}_{\psix \psix} + \eps^{\pseven}_{\psix \psix}
			+ \eps^{\peight}_{\psix \psix}
		\end{array}}
	\end{equation}
\end{example}

\section*{Acknowledgements}
The author is supported by DFG SPP 1388.

\bibliographystyle{abbrv}
\bibliography{refs_def2}

\end{document}